\documentclass[11pt,a4paper]{article}
\usepackage[T1]{fontenc}
\usepackage[utf8]{inputenc}
\usepackage[english]{babel}
\usepackage{lmodern}
\usepackage[margin=1in]{geometry}
\usepackage[numbers,sort&compress]{natbib}
\usepackage{amsmath,amssymb,amsthm}
\usepackage{tikz}
\usepackage{microtype}
\usepackage{xurl}
\usepackage[hidelinks]{hyperref}
\newtheorem{theorem}{Theorem}[section]
\newtheorem{lemma}[theorem]{Lemma}
\newtheorem{conjecture}[theorem]{Conjecture}

\newtheorem{corollary}[theorem]{Corollary}
\theoremstyle{definition}
\newtheorem{construction}[theorem]{Construction}
\newtheorem{definition}[theorem]{Definition}
\theoremstyle{remark}
\newtheorem{example}[theorem]{Example}

\numberwithin{equation}{section}
\hypersetup{
  pdftitle={Prescribed and simultaneous failures of log-concavity in trees with perfect matchings},
  pdfauthor={Heshun Fan},
  pdfkeywords={Independence polynomial, tree, log-concavity, perfect matching},
  pdfdisplaydoctitle=true
}

\title{Prescribed and simultaneous failures of log-concavity in trees with perfect matchings}
\author{Heshun Fan\\
{\small School of Mathematics and Statistics}\\
{\small Huazhong University of Science and Technology}\\
{\small Wuhan 430074, China}\\
{\small \href{mailto:U202410055@hust.edu.cn}{U202410055@hust.edu.cn}}}
\date{}

\begin{document}
\maketitle

\begin{abstract}
We construct trees with perfect matchings and independence sequences
that fail log-concavity near the independence number $\alpha$.
For every integer $k\geq1$, one family fails at $\alpha-k$ and has
a witness of order $O(k^2\log(k+2))$.
For every integer $t\geq2$, a single tree has simultaneous failures
at $\alpha-1,\alpha-3,\ldots,\alpha-(2t-1)$, with a witness
of order $O(t^3)$. All these trees satisfy $\alpha=N/2$, where $N$
is the number of vertices, and contradict Conjecture~4.2 of Ramos
and Sun. A sufficient condition uses a fixed tree and the numbers
of attached leaves. For fixed $k$ or $t$, explicit bounds give geometric decay of
$a_j^2/(a_{j-1}a_{j+1})$ at these indices as the trees grow;
$a_j$ counts independent $j$-sets. The proofs count
independent sets by empty matching pairs. The Lean proof assistant checks formal proofs
of the parameterized results and their quantitative bounds;
exact 90- and 180-vertex examples illustrate the constructions.
\end{abstract}

\begingroup\small
\par
\endgroup

\setcounter{tocdepth}{2}
\tableofcontents

\section{Introduction}\label{sec:introduction}

For a finite simple graph $G$, let $a_j(G)$ be the number of independent
vertex sets of cardinality $j$, with $a_0(G)=1$. Its independence
polynomial is $I_G(x)=\sum_{j=0}^{\alpha(G)}a_j(G)x^j$, where
$\alpha(G)$ is the maximum size of an independent set. Write $\nu(G)$
for the maximum size of a matching. The independence sequence is
log-concave at an index $1\leq j<\alpha(G)$ if
$a_j(G)^2\geq a_{j-1}(G)a_{j+1}(G)$; a failure means the strict
reverse inequality.

Ramos and Sun formulate the following proposed relation between the
position of a failure and the independence number.
For a tree of order $N$, their convention in Section~1.2 replaces
$N/2$ by its floor when $N$ is odd. The conjecture does not separately
specify a range for $\beta$. We use even $N$ and positive integer
$\beta$, with all three indices in the implication below lying in
the support of the independence polynomial.

\begin{conjecture}[{Ramos and Sun~\cite[Conjecture~4.2]{RamosSun2025}}]\label{conj:ramos-sun}
Let $T$ be a tree of order $N$, and put $h=\lfloor N/2\rfloor$.
The proposed implication is
\begin{equation}\label{eq:conjecture}
 a_{h-\beta}(T)^2<a_{h-\beta-1}(T)a_{h-\beta+1}(T)
 \quad\Longrightarrow\quad \alpha(T)\geq h+\beta+1.
\end{equation}
\end{conjecture}

We construct counterexamples with a perfect matching at every prescribed
positive distance from the independence number. The following theorem
also bounds their size and the extent of the failure.

\begin{theorem}\label{thm:main}
For every integer $k\geq1$, set
\begin{equation}\label{eq:parameters}
 q=k+1,\qquad L=\lceil\log_2(k+2)\rceil,\qquad p_0=32qL.
\end{equation}
There is an explicit family of trees $T(q,p)$, indexed by integers
$p\geq p_0$, with a perfect matching and
\begin{equation}\label{eq:order}
 |V(T(q,p))|=2m,\qquad
 \alpha(T(q,p))=\nu(T(q,p))=m=(3q+1)p+q+1.
\end{equation}
Writing $a_j=a_j(T(q,p))$, the index $m-k$ satisfies $1\leq m-k<m$ and
\begin{equation}\label{eq:general-ratio}
 0<\frac{a_{m-k}^{2}}{a_{m-k-1}a_{m-k+1}}
 \leq 2^{-16qL}\left(\frac23\right)^{p-p_0}<1.
\end{equation}
Thus the ratio tends to zero as $p\to\infty$ for each fixed $k$.
The initial tree, obtained by taking $p=p_0$, has order $N_0$ satisfying
\begin{equation}\label{eq:size}
 N_0\leq258(k+1)^2\lceil\log_2(k+2)\rceil
       \leq2064k^2\log_2(k+2).
\end{equation}
\end{theorem}

Taking $h=m$ and $\beta=k$ in~\eqref{eq:conjecture}, its conclusion
would require $\alpha(T(q,p))\geq m+k+1$, whereas
\eqref{eq:order} gives $\alpha(T(q,p))=m$.
For $k=1$, Theorem~\ref{thm:sharp} below proves the sharper threshold
$p\geq6$, and Corollary~\ref{cor:ratio} gives the ratio bound
$49\cdot2^{-p}$. Taking $p=6$ gives the 90-vertex tree.
The general threshold in~\eqref{eq:parameters} is
$128$ in this case.

For simultaneous failures we use a second family, denoted by
$\mathcal T(t,p)$, with any fixed number of odd-distance failures
in a single tree. Its size and the common parameter threshold are explicit.

\begin{theorem}\label{thm:simultaneous}
For every integer $t\geq2$, set
\begin{equation}\label{eq:sim-threshold}
 \ell_t=\lceil\log_2(10t)\rceil,\qquad
 P_t=64(5t+\ell_t)^2.
\end{equation}
For every integer $p\geq P_t$, there is an explicit tree
$\mathcal T(t,p)$ with a perfect matching and
$N=2m$, $\alpha=\nu=m=t(3+7p)$. In this same tree, all $t$ indices
$j_s=m-(2s+1)$, $0\leq s<t$, are distinct and satisfy $1\leq j_s<m$ and
\begin{equation}\label{eq:sim-failure}
 a_{j_s}(\mathcal T(t,p))^2
 <a_{j_s-1}(\mathcal T(t,p))a_{j_s+1}(\mathcal T(t,p)).
\end{equation}
For each fixed $t$, the maximum over $0\leq s<t$ of the positive ratios
$a_{j_s}^2/(a_{j_s-1}a_{j_s+1})$, where $a_j=a_j(\mathcal T(t,p))$,
tends to zero as $p\to\infty$.
At $p=P_t$, the initial tree has at most $6t+201600t^3$ vertices.
\end{theorem}

Theorem~\ref{thm:simultaneous} prescribes the odd offsets
$1,3,\ldots,2t-1$ in one graph. Theorem~\ref{thm:main} allows an
arbitrary single offset and gives the smaller size bound~\eqref{eq:size}.
Both constructions retain a perfect matching and $\alpha=N/2$.

To explain the proof, a perfect matching of size $m$ lets us represent
an independent set of size $m-r$ by its $r$ empty pairs and a choice
of one endpoint in every remaining pair. The edges between pairs
restrict these endpoint choices. The construction of $T(q,p)$ arranges the
restrictions around $q+1$ central pairs. Leaving selected central
pairs empty frees many endpoint choices in their attached branches.
An upper bound for $a_{m-q+1}$ and single positive contributions to
$a_{m-q}$ and $a_{m-q+2}$ then give the required comparison.
In the general construction, central matching pairs are arranged along
a fixed tree, and branch numbers may differ. Unrestricted independent
endpoint choices at $p$ pairs contribute a factor $2^p$.
For each bound on the number of empty central pairs, we maximize
the exponent in these contributions. Theorem~\ref{thm:structural}
shows that a positive second difference of the maxima forces failure
beyond an explicit threshold. For Theorem~\ref{thm:simultaneous}, we
connect several arrangements of three central pairs and, for each
bound, exhibit endpoint choices that attain the corresponding maximum
while respecting the connecting edges.

\subsection{Relation to earlier constructions}\label{sec:related}

Failures close to the independence number and multiple failures in one
tree occur in earlier constructions. We compare their parameters with
the additional requirement of a perfect matching. For a failure at
index $j$, substitution of $\beta=h-j$ in~\eqref{eq:conjecture} gives
the proposed conclusion
\begin{equation}\label{eq:comparison-test}
 \alpha+j\geq2\lfloor N/2\rfloor+1.
\end{equation}
This substitution may give a negative $\beta$ when $j>h$.

Kadrawi and Levit~\cite[Sections~3--4]{KadrawiLevit2023} give infinite
families with $\alpha=N/2+1$ and a failure at $j=N/2$; these satisfy
\eqref{eq:comparison-test}. Galvin~\cite[Theorem~2.1]{Galvin2026}
constructs trees with
$N=1+s+2st$, $\alpha=s(t+1)$, and a failure at $j=st+2$, for
$t\leq s\leq2^{t/16}$ and sufficiently large $t$.
Here $\alpha-N/2=(s-1)/2>0$ and $\alpha+j=N+1$.

Bautista-Ramos~\cite[Definition~1 and Theorem~1]{BautistaRamos2025}
constructs, for each fixed integer $s\geq1$ and sufficiently large
$t$, trees with
\[
 N=2+s(4+6t),\qquad \alpha=1+3s(t+1),\qquad
 j=\alpha-(2r+1)\quad(0\leq r\leq s-1)
\]
at which log-concavity fails. In this family $\alpha-N/2=s$ and
$\alpha+j-N-1=2(s-r-1)\geq0$.
For fixed $s$ and $0\leq r\leq s-1$, the proof of Theorem~1
in~\cite[p.~6]{BautistaRamos2025} gives, at $j=\alpha-(2r+1)$,
\[
 a_j^2=\Theta(2^{(6r+2)t}),\qquad
 a_{j-1}a_{j+1}=\Theta(2^{(6r+3)t}).
\]
The implied positive constants may depend on $s,r$.
Thus the quotient is $\Theta(2^{-t})$ and tends to zero as $t\to\infty$.
Bautista-Ramos, Guill\'en-Galv\'an, and G\'omez-Salgado~\cite[Corollary~4]{BGG2026}
give $(N,\alpha,j)=(19s+6,10s+3,10s+2)$ for $s\geq9$;
again $\alpha>N/2$ and~\eqref{eq:comparison-test} holds.
Their Corollaries~5 and~6 also give two consecutive failures at
$\alpha-1,\alpha-2$ for $(N,\alpha)=(19s+8,10s+4)$, $s\geq16$,
and three at $\alpha-1,\alpha-2,\alpha-3$ for
$(N,\alpha)=(23s+16,12s+8)$, $s\geq13$, respectively.
Counting vertices in their Definition~1 gives these orders;
in both families $\alpha-N/2=s/2>0$, which excludes a perfect matching.
A further construction of Bautista-Ramos~\cite[Definition~1 and Corollary~3]{BautistaRamos2026Spherical}
has arbitrarily many failures in spherically symmetric trees of order
$2^\ell(n+2)-1$, by varying $\ell\geq3$ and then taking odd $n$
sufficiently large.
Their odd order excludes a perfect matching.
The proofs in~\cite[Lemma~3]{BautistaRamos2025} and
\cite[Section~2, Theorem~3 and Lemma~18]{BautistaRamos2026Spherical}
compare dominant coefficient growth rates to force failures; this
comparison is also the basis of our sufficient condition.
Theorems~\ref{thm:main} and~\ref{thm:simultaneous} combine a perfect
matching and $\alpha=N/2$ with, respectively, an arbitrary single
offset or the simultaneous odd offsets, and with explicit parameter
and initial-size bounds. Theorem~\ref{thm:structural} obtains its
sufficient condition from finite core data.

Fang, Lu, Nevo, Yao, and Zheng~\cite[Theorems~1.1--1.2]{Fang2026}
announce unimodality for sufficiently large forests and strict
log-concavity for $N/5\leq j\leq17\alpha/25$.
For $T(q,p)$, the index $j=m-k$ is above that range: for $q=k+1$ and $p\geq1$,
$m\geq4k+6$, so $m-k>17m/25$.
For $\mathcal T(t,p)$ with $p\geq1$, $m\geq10t$ and
$j_s\geq m-(2t-1)>4m/5>17m/25$.
A local failure of log-concavity does not imply nonunimodality.

The representation by matching pairs and oriented edges occurs in the
public implementation~\cite{MatchingLift}. Public Lean
source~\cite{LeanWitnesses} proves that rooted-tree polynomial
calculations count independent sets and provides non-log-concave
witnesses of orders $26$, $28$, and $29$. Its witnesses of orders $26$
and $28$ have $(N,\alpha,j)=(26,14,13)$ and $(28,15,14)$,
respectively, and satisfy~\eqref{eq:comparison-test}.
We use the same general matching-pair representation, with the
construction and all counting arguments specified below.

\noindent\textbf{Organization.}
Section~\ref{sec:preliminaries} constructs the trees and counts independent
sets by empty matching pairs. Section~\ref{sec:general-proof} proves
Theorem~\ref{thm:main} from estimates for these counts.
Section~\ref{sec:sharp} sharpens the estimates when $k=1$, and
Section~\ref{sec:instance} computes the 90-vertex example.
Section~\ref{sec:structural} turns the same counting representation
into a general sufficient condition and applies it to unequal
attachments. Section~\ref{sec:simultaneous} computes the connected
core maxima needed for Theorem~\ref{thm:simultaneous}.
Section~\ref{sec:lean} describes the formalized statements.
Appendix~\ref{app:verification} derives coefficient recurrences on
the original tree and on matching labels, including an explicit
polynomial for the connected 180-vertex example.

\noindent\textbf{Conventions and notation.}
All graphs are finite and simple. The order of a graph is its number
of vertices. A matching is a set of pairwise vertex-disjoint edges;
it is perfect if it covers every vertex. We put $a_j(G)=0$ for
$j>\alpha(G)$ and $I_\varnothing(x)=1$. For a vertex set $Z$, $G-Z$
is the induced graph on $V(G)\setminus Z$, and $N_G[v]$ is the closed
neighborhood of $v$. All graph parameters below are integers;
$\mathbb N=\{0,1,2,\ldots\}$ and $\log_2$ denotes the real logarithm to base $2$.

\section{Preliminaries: the trees and their independent sets}\label{sec:preliminaries}

We first construct the family in Theorem~\ref{thm:main} for the larger
range $q\geq2$, $p\geq1$. A directed tree records the edges between
matching pairs, so the same object describes both the graph and its
independent-set counts.

\begin{construction}\label{con:tree}
For integers $q\geq2$ and $p\geq1$, let $H=H(q,p)$ have distinct
vertices
\[
 z,\quad u_i\ (1\leq i\leq q),\quad
 o_{i,j}\ (1\leq j\leq p),\quad
 b_{i,j}\ (1\leq j\leq2p),\quad d_j\ (1\leq j\leq p),
\]
where both doubly indexed families also have $1\leq i\leq q$.
Its arrows, and no others, are
\[
 u_i\longrightarrow z,\qquad
 u_i\longrightarrow o_{i,j},\qquad
 b_{i,j}\longrightarrow u_i,\qquad
 z\longrightarrow d_j.
\]
Replace each label $v\in V(H)$ by two vertices $U_v,V_v$ and add
the edge $U_vV_v$. For each arrow $v\to w$, add the edge $U_vV_w$.
The resulting undirected graph is $T(q,p)$.
\end{construction}

The underlying undirected graph of $H$ is a star on $z,u_1,\ldots,u_q$
with leaves attached. It is a tree with
$m=(3q+1)p+q+1$ vertices and $m-1$ arrows. The graph $T(q,p)$ is
connected, since each pair is connected and the arrows connect all
pairs. It has $2m$ vertices and $2m-1$ edges, so it is a tree.
The edges $U_vV_v$ form a perfect matching of size $m$. Every
independent set meets each pair in at most one vertex, while
$\{U_v:v\in V(H)\}$ is independent. Consequently
$\alpha(T(q,p))=\nu(T(q,p))=m$. Subsets of this independent set
show that $a_j(T(q,p))>0$ for every $0\leq j\leq m$.

\subsection{Assignments and induced deletion}\label{sec:assignments}

We encode the choice of $U_v$ by $+$, the choice of $V_v$ by $-$,
and an empty matching pair by $0$. An arrow $v\to w$ forbids exactly
the pair of states $(+,-)$. Once the zero labels have been removed,
only assignments of $+$ and $-$ remain to be counted.

\begin{definition}\label{def:assignments}
For a finite directed graph $K$, let $\mathcal C(K)$ be the set of
assignments $\sigma:V(K)\to\{+,-\}$ such that no arrow $v\to w$
has $(\sigma(v),\sigma(w))=(+,-)$. Put
$\Omega(K)=|\mathcal C(K)|$, with $\Omega(\varnothing)=1$.
For $Z\subseteq V(K)$, the graph $K-Z$ retains exactly the arrows
whose two endpoints remain; no constraint passes through a removed
vertex.
\end{definition}

For example, the path with vertices $V_1,U_1,V_2,U_2$ and matching
edges $U_1V_1,U_2V_2$ has the single label arrow $1\to2$.
The three assignments with no empty pair are $(+,+),(-,+),(-,-)$.
For either choice of the empty label, the remaining label has two states;
two empty pairs give one choice. Thus $I_{P_4}(x)=1+4x+3x^2$.

We will group independent sets by their empty pairs. In the lemma,
the original leaves mean the labels $o_{i,j},b_{i,j},d_j$ in
Construction~\ref{con:tree}; the core is $\{u_1,\ldots,u_q,z\}$.

\begin{lemma}\label{lem:deletion}
For $T(q,p)$ and every integer $0\leq r\leq m$, write
$S_r=a_{m-r}(T(q,p))$. Then
\begin{equation}\label{eq:deletion}
 S_r=\sum_{\substack{Z\subseteq V(H)\\|Z|=r}}\Omega(H-Z).
\end{equation}
If $C$ is a subset of the core and $D$ is a set of original leaves,
then
\begin{equation}\label{eq:leaf-extension}
 \Omega(H-C-D)\leq\Omega(H-C).
\end{equation}
\end{lemma}

\begin{proof}
An independent set chooses at most one endpoint in each matching pair.
Let $Z$ be the set of empty pairs. On each remaining pair assign $+$
or $-$ according to the chosen endpoint. The nonmatching edge $U_vV_w$
forbids exactly $(+,-)$ on the arrow $v\to w$. Thus the assignment
lies in $\mathcal C(H-Z)$, and the independent set has size $m-|Z|$.
Conversely, an assignment in $\mathcal C(H-Z)$ selects the indicated
endpoint outside $Z$ and no endpoint in $Z$. Every matching edge has
at most one selected endpoint, and every nonmatching edge either meets
an empty pair or avoids the forbidden pattern. The selected set is
independent. The two constructions are inverse, proving~\eqref{eq:deletion}.

Each original leaf has degree at most one after deleting $C$, and no
two original leaves are adjacent. Every assignment on $H-C-D$ therefore
extends to $H-C$: an isolated restored leaf may take either sign;
on an arrow from its neighbor to the leaf choose $+$ at the leaf,
and on an arrow from the leaf to its neighbor choose $-$.
These choices respect every incident arrow and do not affect one
another. Restriction is consequently surjective, proving
\eqref{eq:leaf-extension}.
\end{proof}

\section{Proof of the prescribed-distance theorem}\label{sec:general-proof}

The sum in~\eqref{eq:deletion} reduces the coefficient comparison to
assignment counts after deletion. We compute the counts for core
deletions, use~\eqref{eq:leaf-extension} to control additional leaf
deletions, and then estimate the resulting ratio for all $p\geq p_0$.

\subsection{Core counts and the three coefficients}\label{sec:core-counts}

Set $B=2^p$. When $u_i$ has sign $+$, its $2p$ incoming leaves are
free and its $p$ outgoing leaves are forced, giving weight $B^2$.
When $u_i$ has sign $-$, the weight is $B$. At $z$ the weights are
$1$ and $B$ for signs $+$ and $-$, respectively.
Deleting a vertex $u_i$ frees its $3p$ leaves, giving a factor $B^3$.

If $C$ consists of $t$ vertices among $u_1,\ldots,u_q$ and does not
contain $z$, then
\begin{equation}\label{eq:core-no-z}
 \Omega(H-C)=B^{3t}\bigl((B^2+B)^{q-t}+B^{q-t+1}\bigr)
 \qquad(0\leq t\leq q).
\end{equation}
Indeed, sign $+$ at $z$ permits either sign at each remaining $u_i$;
sign $-$ at $z$ forces all remaining $u_i$ to have sign $-$.
If $C$ instead contains $z$ and $t$ of the $u_i$, deleting $z$ frees
its $p$ leaves and removes every arrow between remaining core vertices.
Thus
\begin{equation}\label{eq:core-with-z}
 \Omega(H-C)=B^{3t+1}(B^2+B)^{q-t}
 \qquad(0\leq t\leq q).
\end{equation}

Equations~\eqref{eq:core-no-z} and~\eqref{eq:core-with-z} give an upper
bound for $S_{q-1}$, while retaining one core deletion for each of
$S_{q-2}$ and $S_q$ gives lower bounds.

\begin{lemma}\label{lem:coefficient-bounds}
For integers $q\geq2$, $p\geq1$, put $m=(3q+1)p+q+1$ and $B=2^p$.
With $S_r=a_{m-r}(T(q,p))$, one has
\begin{align}
 S_r&\leq\binom mr(2^q+1)B^{2q+r}\quad(0\leq r\leq q-1),
       \label{eq:general-upper}\\
 S_{q-2}&\geq B^{3q-2},\qquad S_q\geq B^{3q+1}.
       \label{eq:general-lower}
\end{align}
\end{lemma}

\begin{proof}
First let $C$ be a core deletion with $c=|C|\leq q-1$.
If $z\notin C$, then $t=c$ in~\eqref{eq:core-no-z}, and
\[
 B^{3c}(B^2+B)^{q-c}\leq2^{q-c}B^{2q+c},\qquad
 B^{q+2c+1}\leq B^{2q+c}.
\]
The inequality $B^{q+2c+1}\leq B^{2q+c}$ uses $c\leq q-1$.
If $z\in C$, put $t=c-1$
in~\eqref{eq:core-with-z}; it gives
\[
 B^{3t+1}(B^2+B)^{q-t}\leq2^{q-t}B^{2q+t+1}
 =2^{q-t}B^{2q+c}.
\]
Hence in both cases
\begin{equation}\label{eq:core-upper}
 \Omega(H-C)\leq(2^q+1)B^{2q+c}\qquad(c\leq q-1).
\end{equation}
For any $r$-element deletion $Z$ with $r\leq q-1$, write $Z=C\sqcup D$
with $C$ in the core and $D$ among the original leaves. Equations
\eqref{eq:leaf-extension} and~\eqref{eq:core-upper}, together with
$c\leq r$, bound its summand by $(2^q+1)B^{2q+r}$.
There are $\binom mr$ choices of $Z$, proving~\eqref{eq:general-upper}.

For $S_{q-2}$, retain the single summand obtained by deleting any
$q-2$ of the $u_i$. Equation~\eqref{eq:core-no-z} bounds this summand
below by $B^{3(q-2)}B^4=B^{3q-2}$.
When $q=2$ the deletion set is empty and this is a bound for $S_0$.
For $S_q$, delete all $q$ of the $u_i$; the resulting summand is
$B^{3q}(1+B)\geq B^{3q+1}$ by~\eqref{eq:core-no-z}.
All summands in~\eqref{eq:deletion} are nonnegative, which proves
\eqref{eq:general-lower}.
\end{proof}

For later use define, for $q\geq2$ and $p\geq1$,
\begin{equation}\label{eq:envelope}
 R_q(p)=\frac{a_{m-q+1}(T(q,p))^2}
 {a_{m-q}(T(q,p))a_{m-q+2}(T(q,p))},\qquad
 F_q(p)=\frac{(2^q+1)^2\binom m{q-1}^{\!2}}{2^p}.
\end{equation}
All coefficients in $R_q(p)$ are positive. Applying
Lemma~\ref{lem:coefficient-bounds} gives
\begin{equation}\label{eq:ratio-envelope}
 0<R_q(p)\leq
 \frac{\binom m{q-1}^{\!2}(2^q+1)^2B^{6q-2}}
 {B^{3q-2}B^{3q+1}}=F_q(p).
\end{equation}

\subsection{The threshold, geometric bound, and size}\label{sec:threshold}

It remains to make $F_q(p)<1$ at an explicit parameter and control
all larger parameters. We use the integer definition
$L=\min\{\ell\in\mathbb N:q+1\leq2^\ell\}$.
Since $2^q\geq q+1$ for $q\geq2$, induction or the binomial theorem
gives $2\leq L\leq q$. Also $q\leq2^L$ and $L\leq2^L$.

\begin{proof}[Proof of Theorem~\ref{thm:main}]
Fix $k\geq1$, let $q=k+1$, and take $p_0=32qL$.
Construction~\ref{con:tree} already gives~\eqref{eq:order}.
We first estimate~\eqref{eq:envelope} at $p_0$. With
$m_0=(3q+1)p_0+q+1$, we have
\begin{equation}\label{eq:m-zero}
 m_0\leq129q^2L\leq2^{8+3L}\leq2^{7L}.
\end{equation}
In~\eqref{eq:m-zero}, the bound $m_0\leq129q^2L$ uses
$3q+1\leq4q$ and $q+1\leq q^2L$; the bound
$129q^2L\leq2^{8+3L}$ uses $129<2^8$ and $q,L\leq2^L$;
and $2^{8+3L}\leq2^{7L}$ follows from $L\geq2$.
Moreover $2^q+1\leq2^{q+1}\leq2^{qL}$, so
\[
 (2^q+1)\binom{m_0}{q-1}
 \leq2^{qL}m_0^{q-1}
 \leq2^{qL+7L(q-1)}\leq2^{8qL}.
\]
Here $\binom{m_0}{q-1}\leq m_0^{q-1}$ follows by encoding each subset
as its increasing tuple. Squaring and dividing by $2^{p_0}$ gives
\begin{equation}\label{eq:initial-envelope}
 F_q(p_0)\leq2^{-16qL}<1.
\end{equation}

To control the increase of the binomial factor, put $a=3q+1$,
$b=q+1$, and $r=q-1$, so $m=ap+b$. For $0\leq i\leq r-1$,
$b-i\geq3$, and hence
\begin{equation}\label{eq:binomial-step}
 \frac{\binom{a(p+1)+b}{r}}{\binom{ap+b}{r}}
 =\prod_{i=0}^{r-1}\left(1+\frac{a}{ap+b-i}\right)
 \leq\left(1+\frac1p\right)^r.
\end{equation}
If $n=2r$, then $p\geq p_0\geq4n$. The binomial theorem and
$\binom ni\leq n^i$ yield
\begin{equation}\label{eq:power-step}
 \left(1+\frac1p\right)^n
 \leq\sum_{i=0}^{n}\left(\frac np\right)^i
 \leq\sum_{i=0}^{n}4^{-i}<\frac43.
\end{equation}
Combining~\eqref{eq:power-step} with~\eqref{eq:binomial-step} gives
\begin{equation}\label{eq:envelope-step}
 \frac{F_q(p+1)}{F_q(p)}
 \leq\frac12\left(1+\frac1p\right)^{2r}<\frac23
 \qquad(p\geq p_0).
\end{equation}
Induction using~\eqref{eq:initial-envelope} therefore proves
\begin{equation}\label{eq:envelope-geometric}
 F_q(p)\leq2^{-16qL}\left(\frac23\right)^{p-p_0}
 \qquad(p\geq p_0).
\end{equation}
Equations~\eqref{eq:ratio-envelope} and~\eqref{eq:envelope-geometric}
prove~\eqref{eq:general-ratio} and its limit. The index is interior
because $m-k=(3q+1)p+2\geq1$ and $k\geq1$.

Finally, at $p=p_0$, equation~\eqref{eq:m-zero} gives
$N_0=2m_0\leq258q^2L$. Since $q=k+1\leq2k$ and
$L=\lceil\log_2(k+2)\rceil\leq2\log_2(k+2)$ for $k\geq1$,
we obtain~\eqref{eq:size}. The parameter $L$ depends only on $k$,
so this is an $O(k^2\log(k+2))$ bound for the chosen initial witness.
\end{proof}

\subsection{Explicit parameters for small prescribed distances}\label{sec:parameter-examples}

Theorem~\ref{thm:main} gives a single formula for every prescribed
distance $k$. Substituting small values separates the role of the
distance from that of the scale $p$. Table~\ref{tab:parameters} records
the initial choices in~\eqref{eq:parameters}; the orders are computed
from~\eqref{eq:order}, rather than from the upper bound~\eqref{eq:size}.

\begin{table}[ht]
\centering
\begin{tabular}{rrrrrrr}
\hline
$k$ & $q$ & $L$ & $p_0$ & $m_0$ & $N_0$ & $m_0-k$\\
\hline
1 & 2 & 2 & 128 & 899 & 1798 & 898\\
2 & 3 & 2 & 192 & 1924 & 3848 & 1922\\
3 & 4 & 3 & 384 & 4997 & 9994 & 4994\\
4 & 5 & 3 & 480 & 7686 & 15372 & 7682\\
\hline
\end{tabular}
\caption{Initial parameters from Theorem~\ref{thm:main}. Here
$m_0=(3q+1)p_0+q+1$ and $N_0=2m_0$.}
\label{tab:parameters}
\end{table}

For instance, distance $k=2$ gives $q=3$ and $m=10p+4$.
The three coefficient indices under comparison are $m-3,m-2,m-1$.
Lemma~\ref{lem:coefficient-bounds} gives, with $B=2^p$,
\begin{equation}\label{eq:distance-two-counts}
 S_1\geq B^7,\qquad
 S_2\leq9\binom m2 B^8,\qquad S_3\geq B^{10}.
\end{equation}
Consequently the middle coefficient satisfies
\begin{equation}\label{eq:distance-two-ratio}
 \frac{a_{m-2}^2}{a_{m-3}a_{m-1}}
 =\frac{S_2^2}{S_3S_1}
 \leq81\binom m2^{\!2}2^{-p}
 \leq2^{-96}\left(\frac23\right)^{p-192}
 \quad(p\geq192).
\end{equation}
The last inequality is~\eqref{eq:envelope-geometric} with $q=3,L=2$.
At $p=192$, this is a failure at index $1922$ in a tree of order
$3848$ with independence number $1924$. In~\eqref{eq:conjecture},
$h=1924$ and $\beta=2$ would require independence number at least
$1927$.

These initial parameters are sufficient choices supplied by one
uniform estimate. The first row can be reduced substantially by
retaining more of the exact core counts. The next section performs
that calculation for $k=1$ and reaches $p=6$.

\section{A sharper estimate at distance one}\label{sec:sharp}

For $q=2$ we can count individual core deletions more precisely than
in~\eqref{eq:general-upper}. This gives counterexamples for $p\geq6$,
well below the general threshold $p_0=128$. Write $T_p=T(2,p)$,
and call the two core labels $u,v$ instead of $u_1,u_2$.

\begin{theorem}\label{thm:sharp}
For every integer $p\geq6$, form a tree from six distinct vertices
$a,b,c,d,e,f$ with edges $ab,cd,ef,af,cf$, by attaching $p$ pendant
paths of length two at each of $a,c,e$ and $2p$ such paths at each
of $b,d$. All new vertices are distinct. This tree is $T_p$ and has
$14p+6$ vertices, a perfect matching, and independence and matching
numbers $m=7p+3$. Its coefficients satisfy
\begin{equation}\label{eq:failure}
 a_{m-1}(T_p)^2<a_{m-2}(T_p)a_m(T_p).
\end{equation}
\end{theorem}

\begin{figure}[ht]
\centering
\begin{tikzpicture}[scale=1.1]
\coordinate (a) at (-2,0);
\coordinate (b) at (-2,1.6);
\coordinate (c) at (2,0);
\coordinate (d) at (2,1.6);
\coordinate (e) at (0,-1.2);
\coordinate (f) at (0,0);
\coordinate (a1) at (-2.9,-.3);
\coordinate (a2) at (-3.8,-.3);
\coordinate (b1) at (-2.9,2);
\coordinate (b2) at (-3.8,2);
\coordinate (c1) at (2.9,-.3);
\coordinate (c2) at (3.8,-.3);
\coordinate (d1) at (2.9,2);
\coordinate (d2) at (3.8,2);
\coordinate (e1) at (0,-2);
\coordinate (e2) at (.9,-2);
\draw[black,line width=0.6pt,dashed] (a)--(f)--(c);
\draw[black,line width=0.6pt,dashed] (a)--(a1) (b)--(b1) (c)--(c1) (d)--(d1) (e)--(e1);
\draw[blue,line width=1.2pt] (a)--(b) (c)--(d) (e)--(f);
\draw[blue,line width=1.2pt] (a1)--(a2) (b1)--(b2) (c1)--(c2) (d1)--(d2) (e1)--(e2);
\foreach \v in {a,b,c,d,e,f,a1,a2,b1,b2,c1,c2,d1,d2,e1,e2}
  \fill (\v) circle (1.6pt);
\node[above right] at (a) {$a$};
\node[above right] at (b) {$b$};
\node[above left] at (c) {$c$};
\node[above left] at (d) {$d$};
\node[right] at (e) {$e$};
\node[above] at (f) {$f$};
\node[below] at (-3.3,-.3) {$p$ arms};
\node[above] at (-3.3,2) {$2p$ arms};
\node[below] at (3.3,-.3) {$p$ arms};
\node[above] at (3.3,2) {$2p$ arms};
\node[below] at (.5,-2) {$p$ arms};
\end{tikzpicture}
\caption{The six-vertex core and one representative of each group of
pendant paths. Each displayed arm represents the stated number of
distinct copies. Thick solid blue edges form the perfect matching;
thin dashed black edges are nonmatching. No arm is attached at $f$.}
\label{fig:core}
\end{figure}
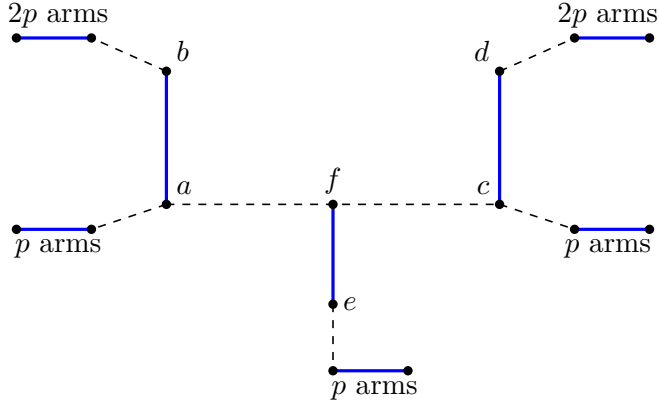

Indeed, identify $(a,b)=(U_u,V_u)$, $(c,d)=(U_v,V_v)$, and
$(e,f)=(U_z,V_z)$. The arrows $u\to z$ and $v\to z$ give $af,cf$;
the attached paths are precisely the leaf pairs of
Construction~\ref{con:tree}. Figure~\ref{fig:core} shows this
description. Thus the graph, matching, and independence number are
already established; we prove the coefficient comparison next.

\subsection{The coefficient comparison}\label{sec:comparison}

We specialize the counts from Section~\ref{sec:core-counts} and sum
single deletions. Put
\[
 w_0=\Omega(H),\qquad w_x=\Omega(H-x),\qquad
 w_{xy}=\Omega(H-\{x,y\})\quad(x\ne y).
\]
Equation~\eqref{eq:deletion} gives
\begin{equation}\label{eq:topthree}
 a_m=w_0,\qquad a_{m-1}=\sum_{x\in V(H)}w_x,\qquad
 a_{m-2}=\sum_{\{x,y\}\subseteq V(H)}w_{xy},
\end{equation}
where the last sum is over unordered two-element subsets.

Set $B=2^p$. Taking $q=2$ and $t=0$ in~\eqref{eq:core-no-z} gives
\begin{equation}
 w_0=(B^2+B)^2+B^3=B^4+3B^3+B^2.
 \label{eq:wzero}
\end{equation}
Taking $q=2$ and $t=1$ in~\eqref{eq:core-no-z}, and using the
symmetry between $u$ and $v$, gives
\begin{equation}
 w_u=w_v=B^3(2B^2+B)=B^4(2B+1).
 \label{eq:wu}
\end{equation}
For $q=2$, equation~\eqref{eq:core-with-z} with $t=0$ and
equation~\eqref{eq:core-no-z} with $t=2$ give, respectively,
\begin{equation}
 w_z=B(B^2+B)^2=B^3(B+1)^2,\qquad
 w_{uv}=B^6(B+1).
 \label{eq:wzuv}
\end{equation}
For $B\geq4$, one has $w_0\leq2B^4$, since
$3/B+1/B^2\leq13/16<1$. Each of the $7p$ leaf deletions has count
at most $w_0$ by Lemma~\ref{lem:deletion}. Summing the three internal
deletions and the $7p$ leaf deletions gives
\begin{equation}
\begin{aligned}
 a_{m-1}
 &\leq 2B^4(2B+1)+B^3(B+1)^2+14pB^4\\
 &=5B^5+(4+14p)B^4+B^3.
\end{aligned}
\label{eq:upper}
\end{equation}

For every integer $p\geq6$,
\begin{equation}
 4+14p\leq \frac{3}{2}2^p.
 \label{eq:threshold}
\end{equation}
The base case is $88\leq96$. If the inequality holds for $p\geq6$,
then $4+14(p+1)\leq2(4+14p)$, because this last comparison reduces
to $10\leq14p$. This proves the induction step for
\eqref{eq:threshold}. Since $B\geq64$, \eqref{eq:upper} now yields
\begin{equation}
 a_{m-1}\leq B^5\left(\frac{13}{2}+B^{-2}\right)<7B^5.
 \label{eq:upperfinal}
\end{equation}
On the other hand, the nonnegative sum in \eqref{eq:topthree} for $a_{m-2}$ includes
the summand $w_{uv}$. Equations \eqref{eq:wzero} and
\eqref{eq:wzuv} therefore imply
\begin{equation}
 a_m\geq B^4,\qquad
 a_{m-2}\geq B^6(B+1)\geq B^7.
 \label{eq:lower}
\end{equation}
As $49<B$, the desired strict inequality follows:
\[
 a_{m-1}^{2}<49B^{10}<B^{11}\leq a_{m-2}a_m.
\]
This proves Theorem~\ref{thm:sharp}.

\subsection{The ratio of the three highest-degree coefficients}\label{sec:ratio}

The estimates~\eqref{eq:upperfinal} and~\eqref{eq:lower} also quantify the failure of
log-concavity. In particular, the ratio in the next corollary can be
made arbitrarily small within the class of trees with perfect matchings.

\begin{corollary}\label{cor:ratio}
For every integer $p\geq6$, let $m=7p+3$ and $a_k=a_k(T_p)$.
Then
\begin{equation}\label{eq:ratio-bound}
 0<\rho_p:=\frac{a_{m-1}^2}{a_{m-2}a_m}<49\cdot2^{-p}.
\end{equation}
In particular, $\rho_p\to0$ as $p\to\infty$ through the integers.
\end{corollary}

\begin{proof}
The coefficients are positive because $\alpha(T_p)=m$.
With $B=2^p$, equations~\eqref{eq:upperfinal} and~\eqref{eq:lower}
give
\[
 0<\rho_p<\frac{(7B^5)^2}{B^7B^4}=\frac{49}{B}=49\cdot2^{-p}.
\]
The upper bound tends to zero.
\end{proof}

\section{The 90-vertex counterexample}\label{sec:instance}

For $p=6$, Theorem~\ref{thm:sharp} gives $N=90$, $m=45$, and
$j=44$. In this section write $a_k=a_k(T_6)$.
We compute its three highest-degree coefficients explicitly;
Section~\ref{sec:lean} describes the formal verification.

\subsection{Exact coefficients}

Allow each matching pair to be empty, to select its $U$ endpoint,
or to select its $V$ endpoint, and give each selected endpoint weight
$x$. Put
\begin{equation}\label{eq:local-polynomials}
\begin{aligned}
 A&=1+2x,\qquad C=1+x,\\
 W&=A^{18}+xA^{12}C^6+xA^6C^{12},\\
 R&=A^{18}+xA^6C^{12}.
\end{aligned}
\end{equation}
A free leaf pair contributes $A$, whereas a leaf pair with one
endpoint forbidden contributes $C$. For either $u$ or $v$, the
three terms of $W$ correspond, in order, to the empty, $+$, and $-$
states. The polynomial $R$ excludes the $+$ state. If $z$ is empty
or has sign $+$, it imposes no restriction on $u,v$; these two cases
contribute $A^6W^2$ and $xC^6W^2$. If $z$ has sign $-$, the $+$
state is forbidden at both $u$ and $v$, giving $xA^6R^2$. Hence
\begin{equation}\label{eq:polynomial90}
 I_{T_6}(x)=(A^6+xC^6)W^2+xA^6R^2.
\end{equation}
Extracting the coefficients in \eqref{eq:polynomial90} gives
\begin{equation}
\begin{aligned}
 a_{43}&=4\,761\,097\,500\,672,\\
 a_{44}&=5\,965\,496\,320,\\
 a_{45}&=17\,567\,744,
\end{aligned}
\label{eq:instance}
\end{equation}
and hence
\begin{equation}\label{eq:instance-difference}
 a_{44}^{2}-a_{43}a_{45}=-48\,054\,595\,706\,911\,981\,568<0.
\end{equation}
The independence number is $45$ by Theorem~\ref{thm:sharp}, while~\eqref{eq:conjecture} at
$N=90$, $h=45$, and $\beta=1$ requires at least $47$.

\section{A structural sufficient condition}\label{sec:structural}

The estimates for $T(q,p)$ depend on which core signs free the attached
leaves. We now allow a general oriented core tree and unequal numbers
of incoming and outgoing leaves. A finite maximum over core states
will control the exponential growth of the actual coefficients.

\begin{construction}\label{con:weighted-core}
Let $K$ be a finite nonempty tree with each edge oriented in exactly
one direction, and put $d=|V(K)|$. Fix nonnegative integers
$\lambda_v,\mu_v$ for $v\in V(K)$, with
$W=\sum_v(\lambda_v+\mu_v)>0$.
For an integer $p\geq1$, attach $\lambda_v p$ distinct outgoing leaves
$v\to o$ and $\mu_v p$ distinct incoming leaves $b\to v$ at each $v$.
All new labels are distinct. Call the resulting directed tree $H_p$.
For every label $x$, introduce $U_x,V_x$ and the edge $U_xV_x$;
for every arrow $x\to y$, add $U_xV_y$. Denote this graph by $G_p$.
\end{construction}

The label tree has $m=d+pW$ vertices. The connected graph $G_p$ has
$2m$ vertices and $2m-1$ edges, so it is a tree. Its pair edges form
a perfect matching. Every independent set meets each pair at most
once, and all $U$ endpoints form an independent set. Hence
\begin{equation}\label{eq:struct-order}
 |V(G_p)|=2m,\qquad \alpha(G_p)=\nu(G_p)=m=d+pW,
 \qquad a_j(G_p)>0\quad(0\leq j\leq m).
\end{equation}
The bijection in the proof of Lemma~\ref{lem:deletion} uses only the
pair edges and the forbidden arrow pattern. It therefore applies
unchanged to $G_p$: for $S_r(p)=a_{m-r}(G_p)$,
\begin{equation}\label{eq:struct-deletion}
 S_r(p)=\sum_{\substack{Z\subseteq V(H_p)\\|Z|=r}}\Omega(H_p-Z)
 \qquad(0\leq r\leq m).
\end{equation}
Here the original leaves are precisely those attached in
Construction~\ref{con:weighted-core}.

\subsection{Core exponents and coefficient bounds}\label{sec:struct-counts}

Keeping all original leaves while deleting core vertices gives an
exact finite sum of powers of $2^p$. The exponent records the leaf
signs that remain free after a legal core assignment is fixed.

\begin{definition}\label{def:core-exponents}
For $C\subseteq V(K)$ and $\sigma\in\mathcal C(K-C)$, define
\begin{equation}\label{eq:core-exponent}
 e(C,\sigma)=\sum_{v\in C}(\lambda_v+\mu_v)
 +\sum_{\substack{v\notin C\\\sigma(v)=+}}\mu_v
 +\sum_{\substack{v\notin C\\\sigma(v)=-}}\lambda_v.
\end{equation}
For an integer $r\geq0$, put
\begin{equation}\label{eq:core-maximum}
 M_r=\max\{e(C,\sigma): C\subseteq V(K),\ |C|\leq r,
                         \ \sigma\in\mathcal C(K-C)\}.
\end{equation}
\end{definition}

The maximum exists: there are finitely many states, and the empty
deletion with every sign $+$ is legal. It is a nonnegative integer
depending only on $K,\lambda,\mu,r$, not on $p$. Its budget is
\emph{at most} $r$ core deletions; the remaining deletions in
\eqref{eq:struct-deletion} may be original leaves.

If a core vertex is deleted, all its attached leaves are free. At
a retained $+$ vertex, incoming leaves are free and outgoing leaves
are forced $+$; at a retained $-$ vertex the roles are reversed.
No original leaves are adjacent to one another. Thus
\begin{equation}\label{eq:struct-core-sum}
 \Omega(H_p-C)=\sum_{\sigma\in\mathcal C(K-C)}2^{p e(C,\sigma)}.
\end{equation}
The number of summands is at most $2^d$. A maximizing core state
also gives a lower bound, after its deletion set is enlarged to
the required size.

\begin{lemma}\label{lem:struct-counts}
For the data of Construction~\ref{con:weighted-core}, let $r\geq0$
and suppose $pW\geq r$. Then
\begin{equation}\label{eq:struct-sandwich}
 2^{pM_r-r}\leq S_r(p)\leq\binom mr\,2^d\,2^{pM_r}.
\end{equation}
The powers in this inequality are interpreted in the positive reals.
\end{lemma}

\begin{proof}
Write a deletion set $Z$ of size $r$ as $C\sqcup D$, with $C$ in
the core and $D$ among the original leaves. Restoring every outgoing
leaf with sign $+$ and every incoming leaf with sign $-$ extends
each legal assignment on $H_p-C-D$ to one on $H_p-C$.
The restored leaves have at most one neighbor and no edges between
them, so all these choices are compatible, even if their core
neighbor was deleted. Consequently
\[
 \Omega(H_p-Z)\leq\Omega(H_p-C)\leq2^d2^{pM_r},
\]
where the second inequality follows from~\eqref{eq:struct-core-sum}
and $|C|\leq r$. Summing the $\binom mr$ terms in
\eqref{eq:struct-deletion} proves the upper bound.

For the lower bound, choose $(C,\sigma)$ attaining $M_r$.
There are $pW\geq r$ original leaves, so choose a set $D$ of
$r-|C|$ of them. Restriction from assignments on $H_p-C$ to
assignments on $H_p-C-D$ has fibres of size at most $2^{|D|}$:
the restriction and the signs on $D$ determine the assignment.
Therefore~\eqref{eq:struct-core-sum} gives
\[
 2^{pM_r}\leq\Omega(H_p-C)
 \leq2^{|D|}\Omega(H_p-C-D)\leq2^r S_r(p).
\]
The last inequality retains the single summand $Z=C\sqcup D$ in
\eqref{eq:struct-deletion}. Division by $2^r$ proves the lower bound,
including the case $pM_r-r<0$.
\end{proof}

\subsection{A positive second difference}\label{sec:struct-threshold}

When $M_{k-1}+M_{k+1}>2M_k$, the product $S_{k-1}(p)S_{k+1}(p)$
has a larger exponential growth rate than $S_k(p)^2$. The polynomial factor in
\eqref{eq:struct-sandwich} can then be controlled at an explicit scale.

\begin{theorem}\label{thm:structural}
Fix the data of Construction~\ref{con:weighted-core} and an integer
$k\geq1$. Define the integer
$\Delta_k=M_{k-1}+M_{k+1}-2M_k$, and suppose $\Delta_k>0$.
Set
\begin{equation}\label{eq:struct-threshold}
 \ell=\min\{l\in\mathbb N:d+W\leq2^l\},\qquad
 n_k=d+k+\ell+1,\qquad p_*(k)=64n_k^2.
\end{equation}
For every integer $p\geq p_*(k)$, the tree $G_p$ has an interior
index $1\leq m-k<m$ and its actual coefficient ratio satisfies
\begin{equation}\label{eq:struct-geometric}
 0<\mathcal R_k(p):=
 \frac{a_{m-k}(G_p)^2}{a_{m-k-1}(G_p)a_{m-k+1}(G_p)}
 \leq2^{-44n_k^2}\left(\frac23\right)^{p-p_*(k)}<1.
\end{equation}
For fixed $K,\lambda,\mu,k$, one has $\mathcal R_k(p)\to0$
as $p\to\infty$ through the integers.
\end{theorem}

\begin{proof}
First suppose $pW\geq k+1$. Then $m=d+pW\geq k+2$,
so all three coefficients are positive by~\eqref{eq:struct-order}
and $m-k$ is interior. Applying Lemma~\ref{lem:struct-counts} at
$r=k-1,k,k+1$ gives
\begin{equation}\label{eq:struct-binomial}
 0<\mathcal R_k(p)=\frac{S_k(p)^2}{S_{k-1}(p)S_{k+1}(p)}
 \leq2^{2d+2k}\binom{d+pW}{k}^{\!2}2^{-p\Delta_k}.
\end{equation}
Put $A=d+W$ and
\begin{equation}\label{eq:struct-envelope}
 E_k(p)=2^{2d+2k}A^{2k}p^{2k}2^{-p\Delta_k}.
\end{equation}
Since $d+pW\leq Ap$ and $\binom{d+pW}{k}\leq(d+pW)^k$,
the right side of~\eqref{eq:struct-binomial} is at most $E_k(p)$.

Write $n=n_k$ and $p_*=p_*(k)$. The definition of $n$ gives
$d+k\leq n$, $k\leq n$, and $\ell\leq n$. Also
$n\leq2^n$, so $p_*=64n^2\leq2^{6+2n}$. Since $A\leq2^\ell$
and the positive integer $\Delta_k$ is at least one,
\begin{equation}\label{eq:struct-initial}
 \begin{aligned}
 E_k(p_*)&\leq
 2^{2d+2k+2k\ell+2k(6+2n)-64n^2}\\
 &\leq2^{14n+6n^2-64n^2}\leq2^{-44n^2}.
 \end{aligned}
\end{equation}
The last comparison uses $n\geq1$ and $14n\leq14n^2$.
Moreover $p_*\geq8k$ and $p_*\geq k+1$. Since $W\geq1$,
every $p\geq p_*$ satisfies the required $pW\geq k+1$.

For $p\geq p_*$, the binomial estimate in~\eqref{eq:power-step},
with exponent $2k$ and $p\geq8k$, gives
\begin{equation}\label{eq:struct-step}
 \frac{E_k(p+1)}{E_k(p)}
 =\left(1+\frac1p\right)^{2k}2^{-\Delta_k}\leq\frac23.
\end{equation}
Induction from~\eqref{eq:struct-initial} yields
$E_k(p)\leq2^{-44n^2}(2/3)^{p-p_*}$.
Together with~\eqref{eq:struct-binomial}, this proves
\eqref{eq:struct-geometric}. The geometric upper bound tends to zero.
\end{proof}

To check $\Delta_k>0$ in Theorem~\ref{thm:structural}, compute
$M_{k-1},M_k,M_{k+1}$ on at most $3^d$ core states, encoding
a deletion by $0$ and a retained vertex by its sign. Neither these
states nor $M_r$ depend on $p$. Thus the test does not require the
independence polynomial of the whole tree. When $\Delta_k\leq0$,
Theorem~\ref{thm:structural} gives no conclusion about log-concavity.

\subsection{Unequal incoming attachments}\label{sec:asymmetric}

The core maxima can be computed for an entire region of unequal
weights. This extends Construction~\ref{con:tree} while retaining
a direct formula for the sufficient condition.

\begin{corollary}\label{cor:asymmetric}
Let $K$ have arrows $u_i\to z$ for $1\leq i\leq q$, where $q\geq2$.
Give $z$ weights $(\lambda_z,\mu_z)=(1,0)$ and give $u_i$ weights
$(1,b_i)$, where the integers $b_i\geq2$ may differ. Put
$B_0=\sum_i b_i$. Then $d=q+1$, $W=B_0+q+1$, and
\begin{equation}\label{eq:asymmetric-maxima}
 M_r=B_0+r\quad(0\leq r\leq q-1),\qquad M_q=B_0+q+1.
\end{equation}
In particular $\Delta_{q-1}=1$, so Theorem~\ref{thm:structural}
applies to these trees with $k=q-1$.
\end{corollary}

\begin{proof}
Let $C$ be the core deletion set. Suppose $C$ contains $t$ branch
vertices and $z$ is retained.
If $z$ has sign $+$, each remaining branch contributes at most
$b_i$, attained at sign $+$, and the exponent is at most $B_0+t$.
If $z$ has sign $-$, every retained branch is forced $-$,
so the exponent is
$q+1+\sum_{u_i\in C}b_i$. For $t\leq q-1$, the sum on retained
branches is at least $2(q-t)\geq q-t+1$, which bounds this exponent
by $B_0+t$. If $z$ is deleted along with $t$ branches, the maximum
exponent is $B_0+t+1$, attained by giving every retained branch sign $+$.
Thus a deletion of at most $r\leq q-1$ core vertices gives at most
$B_0+r$. Deleting any $r$ branches and retaining all other core
vertices with sign $+$ attains this value.

Every exponent is at most the total weight $W=B_0+q+1$.
Deleting all $q$ branches and retaining $z$ with sign $-$ attains
$W$, proving the remaining equality in~\eqref{eq:asymmetric-maxima}.
Its three values at $q-2,q-1,q$ give $\Delta_{q-1}=1$.
\end{proof}

Setting every $b_i=2$ recovers the graph $T(q,p)$ of
Construction~\ref{con:tree}. The coefficient estimates in
Lemma~\ref{lem:coefficient-bounds} give the smaller threshold~\eqref{eq:parameters} and size bound
\eqref{eq:size}; for $q=2$, Theorem~\ref{thm:sharp} further gives
$p\geq6$. As an unequal example, $q=2$ and $(b_1,b_2)=(2,3)$ give
$d=3$, $W=8$, $(M_0,M_1,M_2)=(5,6,8)$ and $p_*(1)=5184$.
Its two branches have incoming leaf counts $2p$ and $3p$, while
both have $p$ outgoing leaves, so a common rescaling of the old
equal attachment counts cannot produce these data.

\section{Several failures in the same tree}\label{sec:simultaneous}

Theorem~\ref{thm:structural} can be applied at several offsets
without changing its core, weights, or scale. To prove
Theorem~\ref{thm:simultaneous}, we construct a connected core whose
maxima have a positive second difference at every required odd offset.

First, for any fixed data in Construction~\ref{con:weighted-core},
the following consequence specifies the common scale.

\begin{corollary}\label{cor:joint}
Let $J$ be a nonempty finite set of positive integers with
$\Delta_k>0$ for every $k\in J$. Put $k_{\max}=\max J$ and
$P=64(d+k_{\max}+\ell+1)^2$, with $\ell$ as in
\eqref{eq:struct-threshold}. For every integer $p\geq P$, the
single tree $G_p$ fails log-concavity at all the distinct interior
indices $m-k$, $k\in J$. For fixed core, weights and $J$,
$\max_{k\in J}\mathcal R_k(p)\to0$ as $p\to\infty$.
\end{corollary}

\begin{proof}
For each $k\in J$, $p_*(k)\leq P$. Apply
Theorem~\ref{thm:structural} to the same $K,\lambda,\mu,p$ at that
$k$. It gives an interior failure index, and different integers $k$
give different indices $m-k$. Each ratio tends to zero.
For any $\varepsilon>0$, take the maximum of the finitely many
scales beyond which the individual ratios are less than
$\varepsilon$. Beyond this scale their maximum is also less than
$\varepsilon$.
\end{proof}

\subsection{A chain of three-vertex blocks}\label{sec:chain-core}

We join copies of the $q=2$ core by arrows between their centers.
The direction of these arrows will allow the states attaining the
desired maxima to coexist in one connected core.

\begin{construction}\label{con:simultaneous}
For an integer $t\geq1$, take distinct vertices $z_i,u_i,v_i$ for
$0\leq i<t$ and precisely the arrows
\begin{equation}\label{eq:chain-arrows}
 u_i\longrightarrow z_i,\qquad v_i\longrightarrow z_i,
 \qquad z_i\longrightarrow z_{i+1}\quad(0\leq i<t-1),
\end{equation}
where the first two arrow families have $0\leq i<t$.
Give each center $z_i$ weights $(1,0)$ and each branch vertex
$u_i,v_i$ weights $(1,2)$. Let $\mathcal T(t,p)$ be the graph
obtained from Construction~\ref{con:weighted-core} with these data.
\end{construction}

The centers form a path with two leaves at each center. Thus the
core is a tree with $d=3t$ and total weight $W=7t$.
Equation~\eqref{eq:struct-order} gives
\begin{equation}\label{eq:chain-order}
 m=\alpha(\mathcal T(t,p))=\nu(\mathcal T(t,p))=t(3+7p),
 \qquad |V(\mathcal T(t,p))|=2t(3+7p).
\end{equation}
We compute $M_r$ for this core, including the effect of the
connecting arrows in~\eqref{eq:chain-arrows}.

\begin{lemma}\label{lem:chain-max}
For Construction~\ref{con:simultaneous} and every integer $r\geq0$,
\begin{equation}\label{eq:chain-maximum}
 M_r=\min\left(7t,\,4t+\left\lfloor\frac{3r}{2}\right\rfloor\right).
\end{equation}
\end{lemma}

\begin{proof}
In the block $\{z_i,u_i,v_i\}$, let $c_i$ count deleted vertices
and let $e_i$ be its contribution to~\eqref{eq:core-exponent}.
A branch contributes $3,2,1$ in states $0,+,-$, respectively;
a center contributes $1,0,1$. With $c_i=0$, a positive center
allows exponent at most $4$, while a negative center forces both
branches negative and gives $3$. With $c_i=1$, deleting the center
gives at most $1+2+2=5$. Deleting a branch gives at most $3+2=5$
at a positive center, or $3+1+1=5$ at a negative center.
For $c_i=2$ or $3$, the total weight bounds $e_i$ by $7$.
These cases prove
\begin{equation}\label{eq:local-exponent-bound}
 2e_i\leq8+3c_i,\qquad e_i\leq7.
\end{equation}
Any globally legal state with at most $r$ deletions has
$\sum_i c_i\leq r$. Summing~\eqref{eq:local-exponent-bound} gives
$e\leq4t+\lfloor3r/2\rfloor$ and $e\leq7t$.

To attain the bound when $0\leq r\leq2t$, write
$r=2s+\epsilon$ with $\epsilon\in\{0,1\}$. In blocks
$0,\ldots,s-1$, delete both branches and give the center sign $-$.
If $\epsilon=1$, delete $u_s$ and give $v_s,z_s$ sign $+$;
the condition $r\leq2t$ ensures $s<t$ in this case.
Give every vertex in the remaining blocks sign $+$.
The centers have a prefix of negative signs followed by positive
signs, so no forward arrow $z_i\to z_{i+1}$ has pattern $(+,-)$.
The internal arrows are also legal. The state has exactly $r$
deletions and exponent
\begin{equation}\label{eq:chain-attainment}
 7s+5\epsilon+4(t-s-\epsilon)
 =4t+3s+\epsilon=4t+\left\lfloor\frac{3r}{2}\right\rfloor.
\end{equation}
For $r\geq2t$, delete every branch and retain all centers with
sign $-$. This legal state uses only $2t$ deletions and has
exponent $7t$. Both parts of~\eqref{eq:chain-maximum} follow.
\end{proof}

The local maximization in the proof can also be read directly from
Table~\ref{tab:block-states}. The coordinates of each displayed state
are $(z_i,u_i,v_i)$, and $0$ denotes a deleted core vertex.
The table concerns one block before the center-to-center arrows are
imposed; those arrows are handled by the ordered choice of maximizing
states in the proof of Lemma~\ref{lem:chain-max}.

\begin{table}[ht]
\centering
\begin{tabular}{cclc}
\hline
$c_i$ & An attaining state & Contribution & Maximum $e_i$\\
\hline
0 & $(+,+,+)$ & $0+2+2$ & 4\\
1 & $(0,+,+)$ & $1+2+2$ & 5\\
2 & $(-,0,0)$ & $1+3+3$ & 7\\
3 & $(0,0,0)$ & $1+3+3$ & 7\\
\hline
\end{tabular}
\caption{Maximum exponent contribution of one three-vertex block
with exactly $c_i$ deletions.}
\label{tab:block-states}
\end{table}

For an odd deletion budget, the construction in the proof of Lemma~\ref{lem:chain-max}
uses $(+,0,+)$ in the partially deleted block, also of exponent $5$,
so that its center remains positive.

For $0\leq s<t$, equation~\eqref{eq:chain-maximum} gives
\begin{equation}\label{eq:chain-three}
 M_{2s}=4t+3s,\qquad M_{2s+1}=4t+3s+1,\qquad
 M_{2s+2}=4t+3s+3.
\end{equation}
Hence $\Delta_{2s+1}=1$, including $s=t-1$, where the last
maximum is attained at the saturation budget $2t$.

\subsection{The common threshold and quantitative bounds}\label{sec:sim-bounds}

The second differences $\Delta_{2s+1}=1$ obtained from
\eqref{eq:chain-three} give the fixed set of offsets $1,3,\ldots,2t-1$. Their individual sufficient thresholds
increase with the offset, so the largest one works for every position.

\begin{proof}[Proof of Theorem~\ref{thm:simultaneous}]
Use Construction~\ref{con:simultaneous} and let
$J=\{2s+1:0\leq s<t\}$. Since $d+W=10t$, the logarithmic
parameter of~\eqref{eq:struct-threshold} is $\ell_t$.
The largest offset is $2t-1$, so Corollary~\ref{cor:joint} has
common threshold
$64(3t+(2t-1)+\ell_t+1)^2=P_t$.
It applies to the single graph $\mathcal T(t,p)$ at every $p\geq P_t$.
Together with~\eqref{eq:chain-order}, it proves the graph assertions,
the distinct interior failure indices, and the fixed-$t$ maximum limit.

More precisely, put $k_s=2s+1$, $n_s=3t+k_s+\ell_t+1$ and
$P_{t,s}=64n_s^2$. For $a_j=a_j(\mathcal T(t,p))$ and
$j_s=m-k_s$, Theorem~\ref{thm:structural} gives
\begin{equation}\label{eq:sim-ratio}
 0<\frac{a_{j_s}^2}{a_{j_s-1}a_{j_s+1}}
 \leq2^{-44n_s^2}\left(\frac23\right)^{p-P_{t,s}}<1
 \qquad(0\leq s<t,\ p\geq P_t).
\end{equation}
Each position retains its own $n_s$ and $P_{t,s}$ in this estimate.

Finally $10t\leq2^{10t}$ implies $\ell_t\leq10t$. Therefore
\begin{equation}\label{eq:sim-size}
 P_t\leq64(15t)^2=14400t^2,\qquad
 |V(\mathcal T(t,P_t))|\leq6t+201600t^3.
\end{equation}
The second inequality substitutes this bound for $P_t$ into
\eqref{eq:chain-order}; it concerns the selected initial tree.
\end{proof}

For $t=2$, one has $\ell_2=5$ and $P_2=14400$.
Thus every integer $p\geq14400$ gives one tree with
$N=12+28p$, $\alpha=6+14p$, and failures at $14p+5$ and $14p+3$.
An exact smaller instance illustrates the difference between this
sufficient threshold and individual parameter values.

\begin{example}\label{ex:180}
For $t=2,p=6$, Construction~\ref{con:simultaneous} gives a tree
with $180$ vertices, a perfect matching and independence number $90$.
Its coefficients, obtainable from~\eqref{eq:rooted-recurrence}
or the polynomial~\eqref{eq:polynomial180}, are
\begin{equation}\label{eq:180-counts}
 \begin{aligned}
 a_{86}&=23\,068\,730\,163\,346\,239\,159\,009\,280,\\
 a_{87}&=46\,506\,972\,783\,451\,140\,259\,840,\\
 a_{88}&=119\,893\,047\,580\,035\,121\,152,\\
 a_{89}&=171\,579\,163\,310\,817\,280,\\
 a_{90}&=304\,089\,070\,043\,136.
 \end{aligned}
\end{equation}
These give
\begin{equation}\label{eq:180-differences}
 \begin{aligned}
 a_{89}^2-a_{88}a_{90}
 &=-7\,018\,756\,060\,810\,230\,465\,580\,742\,079\,414\,272,\\
 a_{87}^2-a_{86}a_{88}
 &=-602\,881\,845\,604\,396\,923\,223\,080\,277\,392\,185\,814\,375\,464\,960.
 \end{aligned}
\end{equation}
Thus the same graph fails at $89$ and $87$. This is a finite
coefficient calculation; Theorem~\ref{thm:simultaneous} with $t=2$
uses $p\geq14400$. The 90-vertex tree in Section~\ref{sec:instance}
is a different, single-block instance.
\end{example}

\section{Lean formalization}\label{sec:lean}

The parameterized inequalities and the 90-vertex example admit formal proofs in
Lean, a proof assistant whose kernel checks that each proof follows
from its definitions and logical assumptions. Here the statements
concern the actual constructed graphs and their independent vertex
sets, including the correspondence between these sets and the
assignment counts used in the paper.

\subsection{The parameterized construction}\label{sec:lean-general}

The formal graph has vertices consisting of a label of $H(q,p)$ and
one of its two endpoints. For this graph $G$, the coefficient of index
$j$ is defined as the cardinality of the collection of all $j$-element
independent vertex sets. In the source, this collection is expressed
as $\{s:\mathrm{Finset}(V(G))\mid G.\mathrm{IsNIndepSet}\ j\ s\}$,
and its cardinality is taken by \nolinkurl{Nat.card}. The formal
\nolinkurl{indepCount} therefore agrees with $a_j(G)$.
The matching number is defined in the project as the supremum of
the numbers of edges in genuine matching subgraphs. The proof shows
that this number is $m$, attained by the displayed perfect matching.

The final theorem \nolinkurl{fixed_k_counterexample} has the assumptions
$k\geq1$ and $p\geq32(k+1)\lceil\log_2(k+2)\rceil$.
It proves treehood, the perfect matching, order $2m$, independence
and matching numbers $m$, the interior failure index, and the strict
inequality between independent-set counts. The theorem
\nolinkurl{fixed_k_strict_failure} extracts the last assertion.

The quantitative results are covered by the following statements:
\begin{itemize}
\item \nolinkurl{actualRatio_le_envelope} proves $R_q(p)\leq F_q(p)$;
\nolinkurl{actualRatio_geometric} proves the positive ratio and
geometric bound in~\eqref{eq:general-ratio}; and
\nolinkurl{actualRatio_tendsto_zero} proves the limit for fixed $q$.
\item \nolinkurl{actualRatio_eq_indepCount} identifies the formal ratio
with the quotient of actual independent-set counts in
\eqref{eq:envelope}.
\item \nolinkurl{fixed_k_initial_size} proves the initial size bounds,
and \nolinkurl{integer_size_bound_le_log_bound} proves explicitly
the comparison of the two right-hand bounds in~\eqref{eq:size}.
\item \nolinkurl{fixed_k_rs42_alignment} identifies $h=m$ and
$j+k=m$ and proves $\alpha(T(q,p))<m+k+1$, the failure of the
proposed conclusion in~\eqref{eq:conjecture} with $\beta=k$.
\end{itemize}

\subsection{The structural and simultaneous results}\label{sec:lean-new}

The general construction uses the same matching-pair representation,
with a finite core and natural-number attachment functions. Its
\nolinkurl{structuralIndepCount} again counts the collection of
independent vertex sets just described. The functions called
\nolinkurl{a} and \nolinkurl{b} in these sources are the outgoing
and incoming multiplicities $\lambda$ and $\mu$ in
Construction~\ref{con:weighted-core}.

In the \nolinkurl{RS42Structural} namespace,
\nolinkurl{structural_counterexample} proves the graph properties
and strict inequality under the core-tree, orientation, weight,
positive-gap and scale hypotheses of Theorem~\ref{thm:structural}.
The statement \nolinkurl{structural_quantitative_conclusion} gives
\eqref{eq:struct-geometric} and convergence of the actual ratio.
The formal maximum $M_r$ uses $|C|\leq r$, as in~\eqref{eq:core-maximum}.
In the \nolinkurl{Asymmetric} subnamespace, \nolinkurl{three_maxima}
proves the values at $r=q-2,q-1,q$ in~\eqref{eq:asymmetric-maxima}
for all $q\geq2$, $b_i\geq2$. For the graphs in
Corollary~\ref{cor:asymmetric}, \nolinkurl{asymmetric_counterexample}
proves the graph properties and strict failure at $m-(q-1)$ for
$p\geq p_*(q-1)$. The theorem \nolinkurl{ratio_geometric} gives
\eqref{eq:struct-geometric} with $k=q-1$, and
\nolinkurl{ratio_tendsto_zero} proves that this actual ratio tends
to zero as $p\to\infty$ with $q$ and the $b_i$ fixed.

In \nolinkurl{RS42Simultaneous}, \nolinkurl{maximum_formula}
proves~\eqref{eq:chain-maximum}, including the legal states attaining
it. The theorem \nolinkurl{simultaneous_counterexample} assumes
only $t\geq2$ and $p\geq P_t$. It binds one graph before quantifying
over all $s\in\{0,\ldots,t-1\}$, and proves the graph properties,
distinct interior indices, and strict count inequalities of
Theorem~\ref{thm:simultaneous}. The identity
\nolinkurl{real_count_definition} identifies each count with
\nolinkurl{Nat.card} of the corresponding \nolinkurl{IsNIndepSet}
collection in that graph. Theorems \nolinkurl{ratio_geometric} and
\nolinkurl{maxRatio_tendsto_zero} give~\eqref{eq:sim-ratio} and the
limit of the maximum for fixed $t$. The statements
\nolinkurl{threshold_bound} and \nolinkurl{initial_size_bound}
give~\eqref{eq:sim-size}. The theorem
\nolinkurl{two_place_counterexample} specializes the failures to
$t=2$, $p\geq14400$. The smaller example in Example~\ref{ex:180}
is supported by a separate finite integer coefficient calculation.

\subsection{The 90-vertex instance and logical assumptions}\label{sec:lean-fixed}

The separate formalization of $T_6$ represents it as
\nolinkurl{SimpleGraph (Fin 90)}, that is, a simple graph on the
vertex set $\{0,\ldots,89\}$. It proves the tree structure,
a perfect matching, and independence number $45$.
The theorem \nolinkurl{graph90_coefficient_eq_polynomial} identifies
every coefficient of the polynomial in~\eqref{eq:polynomial90}
with the corresponding number of independent vertex sets.
The theorem \nolinkurl{graph90_actual_coefficients} proves the exact
counts in~\eqref{eq:instance}, and \nolinkurl{graph90_counterexample}
proves their strict comparison. The integer difference in
\eqref{eq:instance-difference} follows from those counts.

The formalizations were checked with Lean~4.34.0 and mathlib at
commit \nolinkurl{56a9a88447f391672d7aeea117da03b6f7016b0d}.
Their final statements depend only on subsets of the standard assumptions
\nolinkurl{propext}, \nolinkurl{Classical.choice}, and
\nolinkurl{Quot.sound}: equivalent propositions may be identified,
classical choice is available, and equivalent representatives have
equal images in a quotient. The Lean kernel checks the proof terms
relative to these assumptions. The final results use no additional
axioms or unproved placeholders. The formal sources have not been
publicly released.

\section{Discussion}\label{sec:discussion}

Theorem~\ref{thm:main} locates a failure at any prescribed distance
below $\alpha$ while keeping $\alpha=N/2$ and an explicit perfect
matching. The estimates leave room for smaller witnesses. In
particular, when $k=1$, the difference between $p_0=128$ and the threshold $6$
in Theorem~\ref{thm:sharp} comes from replacing the deletion sum
in~\eqref{eq:deletion} at $r=q-1$ by a uniform bound.
Determining sharper thresholds for larger $q$ and the asymptotic behavior
of $R_q(p)$ would require more precise estimates of the same sum
for $S_{q-1}$. The size bound~\eqref{eq:size}
is sufficient for existence and makes no minimality assertion.

Theorem~\ref{thm:structural} makes the required growth comparison
a finite condition on a chosen core and its attachment weights.
It includes the unequal region in Corollary~\ref{cor:asymmetric}.
For the connected chain, the exact maxima~\eqref{eq:chain-maximum}
give simultaneous failures at $\alpha-1,\alpha-3,\ldots,
\alpha-(2t-1)$, with the explicit bounds of
Theorem~\ref{thm:simultaneous}. These specified odd offsets do not
cover arbitrary prescribed sets of positions or consecutive intervals.
Whether the coefficients at other offsets are log-concave is not
decided by a nonpositive value of $\Delta_k$.

\appendix
\section{Independent coefficient calculations}\label{app:verification}

The matching-pair calculation can be checked by recurrences applied
directly to the tree. We record these recurrences to make the alternative calculation
explicit.

\subsection{Recurrences on the original tree}\label{app:rooted-recurrence}

For any graph $G$ and vertex $v$, partitioning independent sets
according to whether they contain $v$ gives
\begin{equation}\label{eq:vertex-recurrence}
 I_G(x)=I_{G-v}(x)+xI_{G-N_G[v]}(x).
\end{equation}
Independent choices in distinct components give the product rule,
with $I_{\varnothing}(x)=1$.

For a tree, a recurrence using a chosen root avoids repeated vertex
deletions. Let $S$ be rooted at $o$, let $\operatorname{Ch}(v)$ be the
set of children of $v$, and let $S_v$ be the subtree induced by $v$
and its descendants. Denote by $F_v^{(0)}(x)$ and $F_v^{(1)}(x)$ the
generating polynomials for independent sets of $S_v$ excluding and
including $v$, respectively. Then
\begin{equation}\label{eq:rooted-recurrence}
\begin{aligned}
 F_v^{(0)}(x)
 &=\prod_{t\in\operatorname{Ch}(v)}
       \bigl(F_t^{(0)}(x)+F_t^{(1)}(x)\bigr),\\
 F_v^{(1)}(x)
 &=x\prod_{t\in\operatorname{Ch}(v)}F_t^{(0)}(x).
\end{aligned}
\end{equation}
Empty products are $1$, so the initial values at a leaf are
$F_v^{(0)}(x)=1$ and $F_v^{(1)}(x)=x$.
At the root, $I_S(x)=F_o^{(0)}(x)+F_o^{(1)}(x)$.
Indeed, excluding $v$ allows either state at each child, whereas
including $v$ excludes every child. Distinct child subtrees have no
edges between them, so their choices multiply; the factor $x$
records the inclusion of $v$.
All these calculations use integer polynomial arithmetic.

\subsection{A three-state recurrence on matching labels}\label{app:label-recurrence}

The matching-pair representation also gives a recurrence directly on
the label tree. It keeps the empty state visible, which is useful for
computing coefficients close to the independence number. Root the
underlying undirected label tree at a label $o$. For a label $v$ and
a state $s\in\{0,+,-\}$, let $Q_v^s(x)$ count legal assignments on
the descendant labels of $v$, including $v$, with state $s$ at $v$.
Each nonzero state contributes a factor $x$ and the zero state
contributes $1$.

For adjacent labels $v,w$, define $C_{vw}(s,t)$ to be $0$ if their
states are forbidden by the arrow between them, and $1$ otherwise.
Thus an arrow $v\to w$ forbids $(s,t)=(+,-)$, whereas an arrow
$w\to v$ forbids $(s,t)=(-,+)$. Any pair involving $0$ is allowed.
Writing $w(0)=1$ and $w(+)=w(-)=x$, one obtains
\begin{equation}\label{eq:label-recurrence}
 Q_v^s(x)=w(s)\prod_{u\in\operatorname{Ch}(v)}
       \left(\sum_{t\in\{0,+,-\}}C_{vu}(s,t)Q_u^t(x)\right).
\end{equation}
At a leaf the triple is $(Q_v^0,Q_v^+,Q_v^-)=(1,x,x)$.
At the root, summing the three states gives the independence
polynomial of the graph obtained by replacing labels with matching
pairs.

To justify~\eqref{eq:label-recurrence}, fix the state at $v$.
Each child subtree then has only one constraint connecting it to $v$,
and no edge joins different child subtrees. The permitted child
states are selected by $C_{vu}$, and the choices in the subtrees
multiply. The factor $w(s)$ records whether the matching pair at $v$
contributes a selected endpoint. The bijection of
Lemma~\ref{lem:deletion} identifies these assignments with independent
sets, so the power of $x$ is their cardinality.

\subsection{Polynomial calculation for connected blocks}\label{app:block-polynomial}

For the chain in Construction~\ref{con:simultaneous}, the recurrence
\eqref{eq:label-recurrence} can first sum all assignments inside each block and then
sum along the center path. Put $A=1+2x$ and $C=1+x$. For every
integer $p\geq1$, define
\begin{equation}\label{eq:block-branch-polynomials}
 \begin{aligned}
 W_p&=A^{3p}+xA^{2p}C^p+xA^pC^{2p},\\
 R_p&=A^{3p}+xA^pC^{2p}.
 \end{aligned}
\end{equation}
The terms of $W_p$ correspond to states $0,+,-$ at a branch label;
$R_p$ omits the positive state. The $2p$ incoming leaves and $p$
outgoing leaves contribute $A$ when their pair is free and $C$ when
one endpoint is forbidden.

Let $E_p,P_p,N_p$ be the polynomials for one block conditioned on
its center being in states $0,+,-$, respectively. The two branches
are independent once the center state is fixed, and therefore
\begin{equation}\label{eq:block-state-polynomials}
 E_p=A^pW_p^2,\qquad P_p=xC^pW_p^2,\qquad N_p=xA^pR_p^2.
\end{equation}
In particular $I_{T_p}=E_p+P_p+N_p$, which at $p=6$ is
\eqref{eq:polynomial90}. For two blocks, the connecting arrow
$z_0\to z_1$ excludes precisely a positive first center and a
negative second center. Removing that contribution from the product
gives the exact identity
\begin{equation}\label{eq:two-block-polynomial}
 I_{\mathcal T(2,p)}(x)
 =(E_p+P_p+N_p)^2-P_pN_p.
\end{equation}
Substitution of $p=6$ yields
\begin{equation}\label{eq:polynomial180}
 I_{\mathcal T(2,6)}(x)
 =\bigl((A^6+xC^6)W^2+xA^6R^2\bigr)^2
       -x^2A^6C^6W^2R^2,
\end{equation}
where $W,R$ are defined in~\eqref{eq:local-polynomials}.
Thus the coefficients in~\eqref{eq:180-counts} can be obtained from
an explicit polynomial as well as from~\eqref{eq:rooted-recurrence}.

For a longer chain, let $D_i^0,D_i^+,D_i^-$ count assignments on the
first $i$ blocks, conditioned on the last center state, and put
$D_i=D_i^0+D_i^++D_i^-$. For fixed $p$, initialize
$(D_1^0,D_1^+,D_1^-)=(E_p,P_p,N_p)$. Adding one block gives
\begin{equation}\label{eq:chain-polynomial-recurrence}
 \begin{aligned}
 D_{i+1}^0&=E_pD_i,\\
 D_{i+1}^+&=P_pD_i,\\
 D_{i+1}^-&=N_p(D_i^0+D_i^-).
 \end{aligned}
\end{equation}
The only omitted transition is again $(+,-)$, and
$I_{\mathcal T(t,p)}=D_t$ for every $t\geq1$. This recurrence
retains the connecting-edge restrictions at each step, in contrast
to multiplying the polynomials of disconnected blocks.

All formulas in this appendix operate over $\mathbb Z[x]$.
For example, the coefficients needed in~\eqref{eq:polynomial180}
follow by finite convolution and
\begin{equation}\label{eq:binomial-extraction}
 [x^j](1+2x)^a(1+x)^b
 =\sum_{i=\max(0,j-b)}^{\min(a,j)}
       2^i\binom ai\binom b{j-i}
 \qquad(a,b,j\in\mathbb N).
\end{equation}
Here $[x^j]$ denotes the coefficient of $x^j$; an empty sum is zero.
Equation~\eqref{eq:binomial-extraction} follows by multiplying the
two binomial expansions and collecting the terms of total degree
$j$. It describes an exact finite calculation of the stated
examples, while the all-parameter inequalities are proved in
Sections~\ref{sec:general-proof}, \ref{sec:structural}, and
\ref{sec:simultaneous}.

\section*{Funding}
This research was self-funded and received no specific grant from
funding agencies in the public, commercial, or not-for-profit sectors.

\section*{Declaration of competing interest}
The author declares no competing financial interests or personal
relationships that could have appeared to influence the work
reported in this paper.

\section*{Acknowledgements}
The author thanks Huazhong University of Science and Technology and
its School of Mathematics and Statistics for providing a supportive
academic environment. The author is also grateful to Tao Gui, Zhiyuan Wang, Maosong Xiang, and Haozhi Zeng
for guidance, encouragement, and providing a workspace for this research.

\end{document}